\documentclass[12pt]{amsart}

\usepackage{amsmath, amssymb}

\newtheorem{thmm}{Theorem}

\newtheorem{obs}{Observation}

\newtheorem{claim}{Claim}[thmm]

\theoremstyle{remark}
\newtheorem*{rem}{Remark}
\newcommand{\br}{\vspace{4mm} \noindent}
\theoremstyle{definition}

\newcommand{\da}{\downarrow}
\newcommand{\free}{\downarrow}
\newcommand{\forks}{\not\downarrow}

\newcommand{\mi}{~$\mathcal{I}$~}
\newcommand{\ua}{\not\downarrow}
\newcommand{\ud}{\downarrow}

\newtheorem*{lemma1}{Lemma}

\theoremstyle{definition}
\newtheorem*{defn}{Definition}

\begin{document}

\section*{Geometry of forking in simple theories}

\begin{center}
\small{\textsc{Assaf Peretz\footnote{Supported by the Clay Mathematics Institute Liftoff fellowship.}} } \\
\end{center}

\br
\begin{quote}
\noindent
\small{\textsc{Abstract}. We investigate the geometry of forking for U-rank 2 elements in supersimple $\omega$-categorical theories and prove stable forking and some structural properties for such elements. We extend this analysis to the case of U-rank 3 elements.}
\end{quote}

\br

\noindent 
Simple theories were defined and initially investigated in 1980 by Shelah ([S2]). In the 90s Kim (and Pillay),  inspired by the work of Hrushovski in the 80s on finite rank cases which showed the possibility of extending the machinery of forking from the stable case to a more general context, developed the basics of simple theories and showed that forking is well behaved in simple theories. Moreover, they showed that simple theories are exactly those theories where the notion of forking (as originally defined by Shelah) is well behaved (e.g. symmetric).


	Early on in the investigation into simple theories, it seemed that the essential behavior of forking was the same as in stable theories and it was conjectured that every instance of forking in a simple theory is witnessed by a stable formula. This conjecture has been formalized 
in various forms as the stable forking conjecture. The truth of this conjecture would imply, among other things, the possibility of a certain lifting of techniques and even
results from the deeply studied stable case to the simple one. Of the various formalizations which the stable forking conjecture has taken, 
counterexamples are known for many. 
Until now the only case where stable forking was known was for 1-based simple theories with elimination of hyperimaginaries.

	Simple theories differ from stable theories as they allow for the independence property while still not having the strict order property 
(every unstable theory has at least one of these two properties). The triangle-free random graph, a theory which has the independence property
and not the strict order property, can be shown to be non-simple. Thus the property of simplicity forms a dividing line inside theories 
without the strict order property. One could conjecture that this line is defined by the fact that relations witnessing the independence 
property cannot be too intertwined with relations witnessing forking. A certain way to formalize this gives a second motivation for the stable forking conjecture, 
which could be described as stating that if two elements fork, then that forking can be witnessed by a relation without the independence 
property. The right form this formulation should take is not yet completely clear. In the $\omega$-categorical case this formulation could 
state that the relation of forking (which is now first-order definable) does not have the independence property.

	Here we will prove that the U-rank 2 elements of an $\omega$-categorical supersimple theory satisfy stable forking. In fact, we will prove 
that the relation of forking itself cannot have the independence property and hence is stable. We will then show how a generalization of the given proof
can be obtained to prove results for higher U-ranks. Specifically, we will investigate the situation for U-rank 3 elements and show results there.

Our method will be to look at the consequences 
which the independence property has for the geometry of forking. We will show that having the 
independence property has surprising consequences for the possible geometry of forking which will force, given simplicity, the relation of 
forking to be stable.

	This paper is derived from the author's Ph.D. thesis under the supervision of Leo Harrington. I would like also to thank
Thomas Scanlon.	

\subsection*{Conventions and Notation}

L is a possibly many sorted language. T is a complete
first order theory in L. We work inside a monster model, M, which is some very saturated model of our theory.  
By an element $a$ we will mean a (possibly imaginary) element of M. 
We will denote first order formulas in our language, of the form $\phi(x,y)$, as
relations $R(x,y)$ (this is to help flesh out the underlying combinatorial structure).

\begin{defn} 
\begin{enumerate}
\item An \emph{IP-sequence} is an infinite indiscernible sequence

\mi which witnesses the
\emph{independence property} for a given formula $R$ (fixed here throughout); 
i.e., for any finite disjoint subsets $A$ and $B$ of \mi, 
\[M~\models \exists x \left( \bigwedge_{a \in A} R(x,a) \land \bigwedge_{b \in B} \neg R(x,b) \right)\]

Recall that a theory $T$ is said to have the independence property 
if it contains a formula having the independence property.

\item Let \mi be an IP-sequence. We say that $c$ is \emph{generic} for \mi with respect to $R$ if there exist infinite
(disjoint) subsets $A, B$ of  \mi such that $\bigwedge_{a \in A} R(c,a) \land \bigwedge_{b \in B} \neg R(c,b)$, $c$ is not algebraic 
in \mi, $ \forall a \in I(a \notin acl(c,I-\{a\}))$ and \mi=$I_1 \cup I_2$ such that $I_1$ and $I_2$ are indiscernible over $c$.  
\br

We say that $c$ is \emph{generic} for $a$ for some $R$ iff $c$ is generic for some \mi with respect to $R$ such that $a\in$\mi. 
We say $c$ is generic in $C$ for $R$, $C$ a finite set, when $c$ is generic for an IP containing all the elements of $C$. 
(Usually $R$ will be obvious and so omitted.) Note that given a formula $R$ with the independence property, we could always get such a 
configuration, i.e., find an $a$ and $c$ such that $c$ is generic for $a$ with respect to $R$.

\item A formula $\phi(x,y)$ is \emph{stable} if there do not exist 
$\langle {a_i} ; i<\omega \rangle$, $\langle {b_j} ; j < \omega \rangle$ such that
$\phi({a_i}, {b_j})$ if and only if $i \leq j$. It is \emph{unstable} otherwise.

\item We say an instance of forking is stable if it is witnessed by a stable formula. i.e.,  
given $a,b$, if $a \not\da_A b$ then there exists a stable formula $\psi(x,y) \in L(A)$ such
that $\psi(a,b)$ and $\psi(x,b)$ forks over $A$.

We say that a theory T satisfies \emph{stable forking} if every instance of forking is stable. 
(For other definitions see references.)

\end{enumerate} 
\end{defn}

We will note that though the definition of genericity seems complicated it is actually a trivial use of the independence property; 
compactness and the definition of the independence property allows us, given an $R$ with the inependence property, to create such a situation
(in the appropriate types in our case), thus allowing us to prepare immediately 
the generic situation which saves a lot of writing in the long run. We will constantly use the generic
situation, which always exists given the independence property, to reach contradictions.

We will only use a few basic results of forking calculus in simple theories:

\begin{enumerate}
\item[0.]
Forking is i) Symmetric: $a \da_A b$ iff $b \da_A a$, and 
ii) Transitive: Suppose $A \subseteq B \subseteq C$, then $a \da_A C$ iff ( $a \da_B A$ and $a \da_B C$).

\item[1.] Let $a$ fork with $b_i$ for $i<n$ and $\{ b_i\}_{i<n}$ be independent. Then
$U(a/\{b_i\}_{i<n})\leq U(a)-n$. (This is clear using symmetry and transitivity).

\item[2.] The Independence Theorem over Lascar strong types. ([W] 2.5.20))
\emph{If $B \da_A C, tp(b/AB)$ and $tp(c/AC)$ do not fork over $A$, and $Lstp(b/A)=Lstp(c/A)$, then there is an
$a \models Lstp(b/A) \cup tp(b/AB) \cup tp(c/AC)$, with $a \da_A BC$.} 

\item[3.] Type-definability of Lascar strong type. ([W] 2.7.9)

\item[4.] If $a\free_A b$ and $Lstp(a/A)=Lstp(b/A)$ then $a$ and $b$ begin a Morley sequence over
$A$. ([W] 2.7.7).

\item[5.] The Lascar inequality (or equality in the finite U-rank case). ([W] 5.1.6)
\emph{$U(a/bA)+U(b/A) \leq U(ab/A) \leq U(a/bA) \oplus U(b/A)$.}

\item[6.] For every $a$ and $A$, $a \da_B A$ for some $B \subseteq A$,
$|B| \leq |T|$.
\end{enumerate}
\br

\section{U-rank 2 stable forking}

In this section we will prove stable forking for $U$-rank $2$ elements, by analyzing
the consequences which the independence property has for forking in simple theories.
Simplifications of the proof can be achieved, but the following proof provides a better background for generalization; also, Theorem 1 is an interesting theorem
on its own right, as it gives information about theories having the independence property but not the strict order property.

\begin{thmm}
Let T be a simple theory. Suppose $e$ is generic for $a$ with respect to $R$, $R(x,a)$ forks over $A$, $U(a/A)=U(e/A)=2$, and $Lstp(b/Ae)=Lstp(a/Ae)$. Then $a \ua_A b$.
\end{thmm}

\begin{proof}
We will omit $A$ in the proof, but all our calculations and analysis are over $A$. As
$e$ is generic for $a$ with respect to $R$, we have by definition of genericity (i.e. the method of getting generic sequences) an IP sequence $I$ with respect to $R$ such that $e$ is generic for $I$, and generic for each $c$, $c \in$I.   

\begin{claim} 
$\forall c \in I (a \ua c)$
\end{claim}

\begin{proof} 
$a,c \in I$, so there are unboundedly many $d$'s which are generic for them
satisfying $R(d,a) \land R(d,c)$ and U(d)=2 (by the indiscernibility 
properties in the definition of genericity). But then
$d \ua a$ and $d \ua c$, so if $a \da c$ then by Fact 1 above we get $U(d/ac)=0$, hence $d$ is algebraic over $ac$ in contradiction to
the way we chose $d$. (We use here our definition of genericity). So $a \ua c$, as desired. 
\end{proof}

\noindent Now suppose for a contradiction that $a \da b$. 
We look at the two possible cases:

1) Suppose $e \da a$. As $e$ is generic for $I$, and as $R$ witness forking, there are unboundadly many $c$'s in I such that  $e \ua c$;
but as $a \ua c$ we get again by fact 1 that $c$ is algebraic in $ea$ which is a contradiction. 

2) Suppose $e \ua a$. 
Now by the Lascar inequalities (fact 5), $U(eab)=U(ab)+U(e/ab)$.  
As we are assuming $a \da b$, we have that $U(ab)=4$, so $U(b/ae)=U(aeb)-U(ae)=
4+U(e/ab)-U(ae)$. 
As we assumed $a \ua e$, $U(ae)=U(a)+U(e/a) \leq 2+1=3$, so 
$U(b/ae) \geq 1=U(b/e)$ hence $b \da_e a$.

As $Lstp(b/e)=Lstp(a/e)$, and as we obtained $b \da_e a$, by fact 4 above we can continue $ab$ to a Morley sequence over $e$, 
$\{a_i\}_{i<w}$.
 
If $c \da_e a$ then by definition of forking there exists some $g$ such that $tp(ge a_i )=tp(cea)$, for all $i$. By fact 6 $c$
cannot fork with every pair of the Morley sequence, so by indiscernibility, and perhaps changing $g$ we get an $f$ such that
$f \da_e ab$ and $f \equiv_{ea_i} g$ for each $i$.  
But now $f \ua a$ and so $f \ua b$. If $a \da b$ then
$U(f/ab)=0$, but $U(f/e)=1$ as being non-algebraic is type definable and hence is in $tp(c/e)$. So $f \ua_e ab$ which is a
contradiction to our assumption, and so we get a contradiction and $c \ua_e a$. 

Now $c$ was just one of infinitely many elements of $I$, and in particular, one of the infinitely many elements in $I$ which
fork with $e$. But $c \ua e$ and
$c \ua_e a$ means $U(c/ae) \leq U(c)-2=2-2=0$ which entails that $c$ is algebraic over $ae$ for infinitely many $c$'s which is a
contradiction. And so $a \ua b$ and we are done.\end{proof}

\vspace{3mm}

\begin{thmm}
Let T be a supersimple, $\omega$-categorical theory. Then the U-rank 2 elements satisfy stable forking.
Moreover, the formula $\phi(x,y) =x \ua y \land U(x)=2 \land U(y)=2$   is stable.
\end{thmm}

\begin{proof}

Let $a$ fork with $d$ over $A$. There is a finite $B \subset A$ such that $a \forks_B d$;
 also as T is supersimple, we can get $U(a/A)=U(a/B)$ and $U(d/A)=U(d/B)$. So  we can assume
without loss of generality that $A$ is finite. This is the only place we use supersimplicity 
rather than simplicity.

$A$ will be omitted in the proof. Let $R$ be a relation such that $aRd$, $R$ witnesses forking, and $xRy$ proves the types of $x$ and $y$.

We note that as we will be dealing with IP sequences and generic elements, elements could be chosen to not be
algebraic over other elements. We will obtain a contradiction to a POSSIBLE configuration, and
hence get a contradiction to any configuration. We will constantly choose elements so as to be
non-algebraic. (This will be much more noticeable in higher U-ranks). 

\begin{obs}
We may assume $R$ has the independence property (IP). 
\end{obs}

\begin{proof} Every unstable formula either has the independence property or some 
conjunction of instances of
the formula has the strict order property (the construction is given explicitly in [S1] Ch. II, \S 4). Now as T is simple, no formula of $T$ 
has the strict order property, and so if $R$ is unstable, it must be because $R$ has the IP. 
Otherwise we are done as then $R$ is stable and witnesses forking.  \end{proof}

We assume $R$ is such that $aRb$ proves $a$ is generic for $b$.

\begin{defn} We define $S$ to be the following relation: $xSb$ iff $tp(x)=tp(a)$, $tp(b)=tp(d)$,
and $b$ forks with some $c$ such that ($xRc$ and $Lstp(c)=Lstp(b)$), and for all $d$ such that
($xRd$ and $Lstp(c/x)=Lstp(d/x)$) then $b$ forks with $d$.  \end{defn} By $\omega$-categoricity (and fact 3), this is
definable.

\begin{claim}
If $xRb$ then $xSb$.
\end{claim} 

This follows from theorem 1.
We get that $b$ forks with each member of the Lascar strong type over $x$ to which $b$ itself
belongs.

\begin{claim} 
S witnesses forking.
\end{claim}

\begin{proof} Let $\{b_i\}_{i<w}$ be Morley, and suppose there exists an $x$ such that $xS(b_i)$ for
all $i<w$. Then, by definition of $S$, there exists some $c$ such that $xRc$, $Lstp(b_i)=Lstp(c)$ and
$b_i$ forks with $c$.  Now notice that by definition we can take the same c for all $b_i$.  (As the
$b_i$s are Morley, they are all of the same Lstp). (There are only a bounded number of Lascar strong types over
$x$, and so some Lascar strong type repeats for infinitely many $c$'s and so $b_i$s, and thus we can get an infinite
subset of our original Morley sequence, which is then still a Morley sequence for which it is the
same $c$).  But this is a contradiction to simplicity (fact 1). \end{proof}

We note that as we got that forking satisfies the definition of $S$, and that $S$ witnesses forking, $S$ is essentially
forking.

\begin{claim}
$S$ does not have the IP.
\end{claim}

\begin{proof} Let $\{b_j\}$ be a set witnessing the IP for $S$. We can assume the sequence $\{b_j\}$
is indiscernible, and as $tp(b)$ is proved by $S$, that $b_0=b$.  Notice that as the
$b_j$s are indiscernible, all of them are of the same Lascar strong type. There exists an $x, xSb_0 \land \neg xSb_1$. 
As $x\free b$, there is some $c$, $Lstp(c/x)=Lstp(b_0/x)$ and $xRc$ such that $b_1$ does not fork
with $c$. 

\begin{rem} Given an IP set $\{b_i\}_{i<w}$, all of them indiscernible (we can achieve this by
compactness, and as $S$ proves type of $b$, can assume $b_0=b$), then $\forall i,j (b_i \ua b_j)$.  \end{rem}

\begin{proof} Otherwise we have infinitely many $x$'s such that $xSb_0$ and $xSb_1$ and $b_0\free
b_1$. So by fact 4    $U(x/b_0 b_1)=0$, so x has to be algebraic over $b_0$ and $b_1$.
Contradiction. \end{proof}

Now, by the IP and the assumption regarding $\omega$-categoricity it is easy to see that there
must be an $x$ and a $y$ such that $Lstp(x/b_0)=Lstp(y/b_0)$ and $xSb_0$ and not $xSb_1$ and
$ySb_0$ and $ySb_1$. (e.g. look at the $\omega$ many subsets of $\omega$ which are $\omega -j$ for each j. then
there is $x_j$ which witnesses $x_jSb_k$ iff not $k=j$. Now at least 2 of them have the same Lstp
over $b_0$, and then can assume $j=1$ without loss of generality).  For this $x$ there is a $c$ as before. $x$
could also be chosen to have infinitely many $b$'s from the IP of the same Lstp over $x$ as
$b_0$ and hence will also fork with $c$.  Now, by our previous remark, $b_0\forks b_1$; by our choice
of $c$, $b_0\forks c$, and also $c\free b_1$, hence by fact 1 we get $U(b_0/cb_1)\leq 2-2=0$
which is a contradiction. 

\end{proof}

And so we are done proving our Theorem. Notice that if we take $R(x,y):= x \ua y \land U(x)=U(y)=2$ (which is
definable in an $\omega$-categorical theory), then if the relation $R$ has the independence property, an extension of $R$ which also states
the types of $x$ and $y$ for some $x$ and $y$ has the IP, while still witnessing forking and the ranks; but as we proved that cannot happen, the
relation of forking, plus the ranks, cannot have the independence property. 
 
\end{proof}

The fact that we got a precise formula (moreover that it is forking) is very useful when one comes to extend this result to higher U-ranks  
(e.g., for a proof by induction), and allows the result itself to be used and not just the methodology of the proof.

\begin{rem} 
(1) I will mention that we could prove claim 2.3 using the independence theorem over Lascar strong types. Such a proof, though longer in 
this case,
allows for certain generalizations to higher U-ranks by changing our elements to tuples.

(2) Using regular type machinery one could take our U-rank 2 elements to be pairs of U-rank 1 elements.
\end{rem}

\section{The U-rank 3 case of stable forking}

In this section we will prove that forking between elements of U-rank 3 is either stable, or is witnessed by an IP sequence with a particular U-rank configuration.
Our elements in this section are of U-rank 3. In this case we have several possible independence property sequences. Let $x$ be generic to an IP-sequence $I$ with respect to $R$ which witnesses forking, and $a,b,c$ distinct elements in $I$. We will divide the possible IP's into 4 cases:
\begin{enumerate}
\item U(b/a)=1. We name it 3-1.
\item U(b/a)=2. We name it 3-2.
\item U(b/a)=3 and U(d/ab)=2. We name it 3-3-2.
\item U(b/a)=3 and U(d/ab)=1. We name it 3-3-1.
\end{enumerate}

In our proof $R$ will be taken to be the forking relation itself (and stating $R(x,y)$ implies $U(x)=U(y)=3$). We will show that in this case, the only possible IP-sequences is of the form 3-3-1.
We note that the 4 cases defined above are the only possible IP sequences, as if a sequence were of the form 3-3-3 then there would be only a finite number of elements of U-rank 3, forking with any 3 elements from the IP, in contradiction of it being an IP. Likewise we cannot have U-rank 0 in IP sequences due to indiscernibility.   

If we are interested in having stable forking for elements of U-rank $\leq 3$ then it is sufficient to prove it for U-rank 3 elements, as if $U(a)<3$ we look at the following structure. We partition the universe into 2 disjoint unary predicates P and Q. P we take as our original structure, while Q is an infinite set with no relations. There  are no relations between the universe of the 2 predicates. We now add to $a$ the needed number of elements from Q to make it U-rank 3 (it is actually clear that in this case $U(a)=2$ and we will add an element, $q$, from Q to $a$. As the only relation with regard Q is equality, we get $U(aq)=3$. We do this with different $q$'s for each member of the IP-sequence. Now if we prove stable forking in this case, it would translate immediately to stable forking in P. 
We will later show that in the case of both elements having U-rank $\leq 3$ but one of the elements having U-rank 2, we have stable forking.

\begin{thmm} T supersimple, $\omega$-categorical.  Let $U(x)=U(a)=3$ and $x$ forks with $a$. Then either the forking is stable or, with respect to the relation $R(z,y):= z \ua_A y \land U(z)=U(y)=3$, $x$ is generic for $a$ and the only IP-sequences which witness this genericity are 2-independent. Furthermore, the sequences are of the form 3-3-1. \end{thmm}

\begin{proof}
Similarly to the Rank 2 case, without loss of generality we can drop the base set $A$. 

We first look at the case where the IP sequences are not 2-independent (i.e. of form 3-1 and 3-2).

\begin{claim} \label{pairs} Suppose that there is an IP-sequence \mi whose elements are not 2-independent. Then whenever $x \ua a$,
$x \ua b$, all generic, and $Lstp(a/x)=Lstp(b/x)$ then $a \ua b$. \end{claim}

\begin{proof} 
We assume towards contradiction that $a \da b$.
There are two subcases. 

\noindent\emph{Subcase 1}.
Let \mi be of the form 3-1. 

Let $x$ be generic for \mi. Recall that \mi is IP with regard to the relation of forking. Let $x$ fork with $a$ and not fork with
$d$. So $U(xda)=U(xd)+U(a/xd) =7$ (as $a$ is not algebraic over $xd$) = $U(ad) + U(x/ad)=4+U(x/ad)<7$ (as $x$ forks with $a$),
contradiction. So no such \mi exists.

\vspace{3mm} 
\noindent\emph{Subcase 2}. Let \mi be of the form 3-2.

We have $xRa, xRb$ where $a \ud b, Lstp(a/x)=Lstp(b/x)$, as well as an IP \mi of form 3-2 for which $x$ is generic. In \mi, let $xRa \land xRe \land xRg \land xRh \land \neg xRd \land \neg xRf$.  Now $U(a/x)=U(x/a)=2$ as otherwise by genericity (and Lstp(b/x)=Lstp(a/x)) $U(xab)=U(a)+U(x/a)+U(b/xa) \leq 3+1+1=5$ while $U(ab)=6$ contradicting the Lascar inequality. 

Suppose towards contradiction that $U(x/ae)=2$. 

By U-rank calculations this implies $U(e/ax)=2$. So $e$ does not fork with $a$ over $x$.  But then we can extend the type of $e$ over
$x$ to both $a$ and $b$, and by the Independence theorem over Lascar strong types ( U-rank calculations prove $a \da_x b$) we get an $l$, $tp(l/xa)=tp(e/xa) \land tp(l/xb)=tp(e/xb)$, which does not fork with $ab$ over
$x$, but which forks with $a$ and with $b$ (which are independent). But now $l$ forks with $a$ and with $b$, and $a$ is
free from $b$, so $U(l/ab) \leq 1 \land U(l/x)=2$ contradiction. 

So $U(x/ae)=1$ and $U(e/ax)=1$. Notice that we are not
using here the fact that $a,e$ are part of an IP-sequence for which $x$ is generic, but only that $Lstp(a/x)=Lstp(b/x)$. (The case where
$a,e$ are not necessarily part of such an \mi allows also for $U(x/ae)=0$.) So $e \ua_a x$, $U(e/a)=2$ and $U(e/ax)=1$. 

Now, $U(xaeg)=6+U(g/xae)=7$ (for $g$'s not algebraic over $xae$, which means for almost all in the IP).

$U(xaeg)=5+U(g/ae)+U(x/aeg)=6+U(g/ae)$ so $U(g/ae)=1$. (for $x$ not algebraic over $aeg$, hence for all generic $x$ for the IP, but
in particular could always find such an $x$, and as the conclusion does not mention $x$, the conclusion follows).  This is not
related to whether $x$ is connected or not to $a$ or $e$ or $g$. So for every $a,e,g$ in \mi, $ U(e/a)=2$ and $U(g/ae)=1$.

Now $U(g/xdf)=1$, so $x$ does not fork with $g$ (or $h$) over $df$, also $g$ does not fork with $h$ over $df$ (otherwise $h$ would be
algebraic over $gdf$), so by the independence theorem over Lascar strong types, we get an $z, z \ud_{df} gh$, so $U(z/df)=1 $ as $U(z/gh)=1 $
(recall this is the case whether $z$ is connected to an IP-sequence containing them or not, as long as $z$ forks with both of them).
Also notice that $x$ or $df$ could be chosen so that $x$ is not algebraic over $df$. (We constantly talk about the generic situation). 

(Notice this is the same computation for $e$ instead of $f$).

So $U(zdf)=5+U(z/df)=6$ but $U(zdf)=6+U(f/zd)$ and so we obtain a contradiction, as $f$ is not algebraic over $zd$.

\end{proof}

\begin{claim} The relation R cannot have the IP with respect to an IP sequence which is not 2 free. \end{claim}

\begin{proof} 

We prove this similarly to the U-rank 2 case (or rather, the remark after the proof where an alternative, a bit longer, proof is mentioned).
\br

\begin{lemma1} If $x$ does not fork with $b$, then there exists a $c$ in each Lascar strong type over $x$ such that $b$ does not fork with $xc$ (even
infinitely many such $c$'s). 
\end{lemma1} 

\begin{proof}
Suppose not. Then $b$ is free from $x$ but forks with $xc$ for all $c$ of some Lascar strong type over $x$. So $b$ forks with $c$ over $x$ for
all such $c$'s. Take a Morley sequence in $b$ over $x$. As $x$ is free from $b$, then that sequence is actually Morley over the
empty set (or whatever base set we're working over). But now as there are only boundedly many Lascar strong types and the series can be taken to be
arbitrarily long, we can assume there is some Lascar strong type over $x$ such that for every $c$ in this Lascar strong type over $x$, $b_i$ forks with $xc$, for each $b_i$ in the Morley sequence. But then $xc \ua_{{b_i}_{i<\alpha}} b_\alpha$ contradicting simplicity.
\end{proof}

Let \mi be an IP, $x$ generic for \mi and $a,b$ be in $I$. By our previous claim we get $U(b/a)=2$ for each $a,b$ in \mi. Let $xRa \land x(\neg R)b$.
Then as $x$ does not fork with $b$, by our previous lemma there exists some $c$ $Lstp(c/x)=Lstp(a/x)$ such that $b$ does not fork with $xc$. By our first claim 
$a \ua c$. 

We will now use the Independence Theorem over Lascar strong types to get a contradiction.
We first prove the requirements:

0) As before, we can assume $U(x/a)=2$ as otherwise $U(xab)=U(a)+U(x/a)+U(b/xa) \leq 6$ while $U(xb)=6$ so $U(a/xb)=0$ but this is true for
infinitely many $a$'s so we get a contradiction, as we can assume $x$ is connected to infinitely many $a$'s (or at least could always
choose such an $x$).

1) We can get a $y$ such that $y$ does not fork with $b$ over $a$, and $Lstp(y/a)=Lstp(x/a)$.

We take a generic $y$ for \mi which forks with $a$ and $b$. Suppose $y$ forks with $b$ over $a$. As $ \forall d \in I(U(d/a)=2)$ and as $y$ is of U-rank 2 over $a$ we get by our result on U-rank 2 elements that the
forking formula has to be stable. In particular it cannot divide the (indiscernible over $a$) IP into 2 infinite parts ([S1] Ch. II Theorem 2.20 proves that if $\phi(z,t)$ is stable and $J$, some indiscernible sequence, then the set of $j$'s in $J$ such that $\phi(a,j)$ is either finite, or cofinite (in $J$) for all $a$'s) . So there are infinitely
many $d$'s which $y$ does not fork with, but which it forks with over $a$ (or we can find a different $b$ that $x$ forks with and $y$
does not fork with over $a$). But now $U(xda)=U(xd)+U(a/xd)=6+1=7$ (the 1 is as $a$ forks with both $x$ and $d$ which are
independent), while $U(xda)=U(da)+U(x/da)=5+1$ (the 1 is by our assumption of forking over $a$) Contradiction.
The fact that we can get $Lstp(y/a)=Lstp(x/a)$ can now be easily seen by counting (same argument as in the U-rank 2 case).

2) $x$ does not fork with $c$ over $a$:  Suppose it does fork. Then $U(xabc)=U(bxc)+U(a/bxc)\geq 8+1=9$ as $a$ cannot be algebraic (or
at least chosen as not to be algebraic from having the IP). But also $U(xabc)=U(ac)+U(x/ac)+U(b/xac)\leq 5+1+2=8$ ($U(ac)\leq 5$ by
our previous claim that $x$ cannot fork with 2-independent elements of the same Lstp over it).  contradiction.

3) $b$ does not fork with $c$ over $a$: If it does, then $U(abc)<7$ but $U(bc)=6$ so $U(a/bc)=0$ hence $a$ is algebraic over $bc$
which is a contradiction, as again, can choose a configuration where that does not happen.

So we can now use the independence theorem over Lstp and get a $z$ such that $z$ forks with $b$ and with $c$ and such that $z$ does not
fork with $bc$ over $a$. But as $U(z/a)=2$ and $z$ forks both with $b$ and with $c$ and $b$ and $c$ are independent, then
$U(z/bc)<2$ hence $z$ forks with $bc$ over $a$ and we get a contradiction.

And so if there is an IP sequence which is not 2-independent, then the forking relation is stable.

\end{proof}

We now look at the case where the IP sequence is 2-independent. In the 3-3-2 case we can use our U-rank 2 result and not just its methods.

\br
Suppose there exists an IP sequence which is 2-independent and such that for $a,b,d$ in \mi $U(d/ab)=2$ (i.e. of form 3-3-2). Then we look at the
new sequence \mi over $a$. Over $a$ any 2 elements in \mi fork with each other. Now this is still an IP sequence as if $x$ forks
with $b$ then it forks with $b$ over $a$ as $a$ and $b$ are independent. If $x$ is free from $b$ we show it is free from $b$ over
$a$. Let $x$ fork with $a$ and $d$ but not with $b$. Then $U(xabd)=U(abd)+U(x/abd)=8+1=9$ so $9=U(xabd)=U(xb)+U(a/xb)+U(d/xab)$.
Now as $U(xb)=$ we get $U(a/xb)=2$ so $U(xab)=U(a/xb)=8$ hence $U(b/xa)=8-U(xa)=3$ and we get $U(b/xa)=U(b/a)$ so $x$ does not fork
with $b$ over $a$ as we wanted to show. 

So we got that if there exists such an IP sequence with respect to forking, then there exists one which is not 2-independent, and as we
know no such sequence exists, then also the independent sequence cannot exist. 

Hence we obtain that the only possible IP-sequence, with respect to forking, is of the form 3-3-1 as desired.
\end{proof}

For a discussion of the 3-3-1 case, see [P].

We finish this section with the following remark (T is as usual):
\begin{rem}
Suppose $x$ forks with $a$ and $U(x)=2 \land U(a)=3$, then this is an instance of stable forking.

Let $R(z,t):= z \ua t \land U(z) \land tp(t)$.
As explained before theorem 6, we add an element $q$ to $x$ such that $U(xq)=3$. Then $x$ still forks with $aq$. If $R$ is not stable, we have an $xq$ and an \mi for which $xq$ is generic for, such that the elements of \mi are of the type of $a$. As we know the only possible IP is of form 3-3-1, and that $xq$ is generic for some $a,b$ and $a,c$ and $\{a,b,c\}$ is independent and $xq$ forks with all 3 elements. But now as the forking cannot occur with respect to the $q$'s as the only way to fork with them is with equality, $x$ forks with $a,b$ and $c$ which are independent, while $x$ is U-rank 2, which is a contradiction.
And so we are done.
\end{rem}

\section{Generalizing Theorem 1}

In this section we will prove a generalization of Theorem 1. This theorem will show again a consequence of a formula both having the
independence property and witnessing forking in a simple theory. We will prove the theorem for U-rank 3 elements as well as remark on a possibility for generalization
to arbitrary finite U-rank. One of the reasons for an interest in this theorem is that it, again, gives us information on the area of non-simple theories without the 
strict order property. Beside that, this theorem gives us structure information inside simple theories, as even given stable forking, still forking can exist with respect 
to a formula with the independence property (stable forking only states that a stable formula witnessing the forking exists).

\begin{defn} 
We say an IP sequence $I$ is  \emph{sound} if for $a \in I$ : $U(a)=n$ for some $n< \omega$ and every size $n-1$ subset, $J$, of $I$ is independent, and $U(a/J)=n-1$ for $a \notin J$.
\end{defn}

For U-rank 3, an IP sequence is sound if for $a,b,d \in$\mi, 
$U(a)=U(b/a)=3 \land U(c/ab)=2$.

We will prove that, as in the U-rank 2 case, if a formula $xRa$ both witnesses forking and has the IP with respect to a sound IP sequence, then $x$ cannot fork with 3 
independent elements of U-rank 3 of the same Lascar strong type over it. This can be generalized to the general finite U-rank situation by adding a base set, C, and 
then demanding  that $a$, that $b$ and that $c$ continue the base set to a generic IP for $x$. This is essentially the same proof, only with trivial alteration of the 
U-rank calculations.

\begin{thmm} Let T be a simple theory. Let $x \ua a$, where $U(x)=U(a)=3$; let the formula $R$ witness the forking and have the independence property where $x$ generic in $a$ with respect to a sound IP sequence.  Then, if $xRa \land xRb \land xRc$, $a,b$ and $c$ of the same Lascar strong type over $x$, then $a \ua bc$. 
\end{thmm}

\begin{proof}

We divide the proof to 3 subclaims, from which the theorem follows.

\begin{claim} 
\label{vee} 
Suppose $a$ and $b$ are elements of a sound IP-sequence \mi for which
$x$ is generic. suppose $a$ and $c$ be elements of another such sequence, that $Lstp(b/ax)=Lstp(c/ax)$ 
and that $a,b,c$ are pairwise independent. Then $\{a,b,c\}$ is not independent. 
\end{claim}

\begin{proof}

Suppose for a contradiction that $\{a,b,c\}$ are independent in the sense of forking. Let $d$ continue
$a,b$ to a sound IP-sequence, $I$.
We remember that for any $d \in $\mi, $U(d/ab) = 2$ and $U(b/a)=3$.

\vspace{2mm} Suppose that $d \ud_ax b$. Now as $Lstp(b/ax)=Lstp(c/ax)$, we can find a Morley sequence $\{ b_j \}$ over $ax$ whose first
two elements are $b,c$. By indiscernibility, for each member of this Morley sequence we get a corresponding $d$; in particular we
can find such a $d_i, d_j$ for some $b_i, b_j$ so that $Lstp(d_i/ax)=Lstp(d_j/ax)$. By boundedness of the number of distinct Lascar strong types and 
type-definability of the independence property, we may assume $i=0, j=1$, i.e., these are the corresponding $d$'s for $b$ and $c$ respectively. 

By assumption, $d_0 \ud_ax b$ and so also $d_1 \ud_ax c$, as $Lstp(cd_0/ax)=Lstp(cd_1/ax)$. By $U$-rank 
calculations,  $b\ud_ax c$ as:  $U(xabc)=U(abc)+U(x/abc)$ (by the Lascar inequality for the finite U-rank case)
$= 3^2+0$ ($x$ has to be algebraic over $abc$ as they are independent and $x$ forks with each of them) $=U(x)+U(a/x)+U(c/xa)+U(b/xac)$;
so as $U(x)=3, U(a/x) \leq 2$ and the last 2 terms are each $\leq 2$, hence $U(b/xac)$ has to equal $2$ so $U(b/xac)=2=U(b/xa)$ hence $b \ud_{xa} c$.

By the independence theorem over Lascar strong types, we can find an $e$ such that $e \ud_{ax} bc$ 
and $Lstp(e/xa)=Lstp(d_0/xa)$ and $(e/xab)=(d_0/xab)$ and $(e/xac)=(d_1/xac)$. As $bc \ud_x a$ we have
$e \ud_x abc$. But now in particular $e \ua ab$ and $e \ua ac$, and as $b \ud_a c$ we get $e \ua_ab c$
(for otherwise $e \ud_ab c \land c\ud ab \rightarrow e \ud_a b \land a \ud b \rightarrow e \ud ab$, a
contradiction). So $U(e/abc) \leq 1$, while $U(e/x) = 2$, which gives $e \ua_x abc$, a contradiction.

\vspace{2mm}
So it must be the case that $d \ua_ax b$. Now, $U(xabd)=U(x)+U(a/x)+U(b/ax)+U(d/xab)=3+2+2+1=8=U(abd)+U(x/abd)=U(ab)+U(d/ab)+U(x/abd)=6+2+U(x/abd)$, so $U(x/abd)=0$, meaning that $x$ 
is algebraic over $abd$.
But $a,b,d$ are part of an IP-sequence generic for $x$ ($a$ here refers to the set of $a_i$'s which form part of this sequence), which means that $x$ is not algebraic with the IP sequence, so we get a contradiction. \end{proof}

\begin{claim}
\label{triangle}
Assume $x$ is generic for $a$ with respect to $R$ with a sound IP sequence, $xRa$ and that $Lstp(c/ax)=Lstp(b/ax) \land Lstp(b/x)=Lstp(a/x)$; furthermore that $\{a,b,c\}$ are pairwise independent. Then $a \ua bc$.
\end{claim}

\begin{proof}
As $Lstp(a/x)=Lstp(b/x)$, we can continue $a,b$ to a Morley sequence $\mathcal{J}$ over $x$. $x$ is generic in $a$ with respect to a sound IP sequence \mi. Let $e \in$\mi, then $e \ud a$. 
By rank calculations again, $e \ud_x a$ so we can find such an
$e$ simultaneously for all the elements of any Morley sequence over $x$, in particular for $\mathcal{J}$. Now as there are
unboundedly many elements in the sequence, we can find 2 such that $Lstp(a_i/ex)=Lstp(a_j/ex)$, and since the
$a_k$'s are Morley over $x$, $e$ cannot fork with all such pairs over $x$; so there exist a pair $a_i, a_j$ so 
that $e\ud_x a_i a_j$. But then by indiscernibility and type-definability of our assumptions, we can take $a_i$
to be $a$ and $a_j$ to be $b$. But now $e,a,b$ satisfy the hypotheses of Claim \ref{vee}, playing the parts
of $a,b,c$ respectively. So by Claim \ref{vee}, $e \ua ab$, and by our choice of $e$, $e \ud_x ab$. 

As in Claim \ref{vee}, suppose $e \ud_{ax} b$. Then the Independence Theorem over Lascar strong types gives a $d$
such that $d \ud_{ax} bc$, $Lstp(d/abx)=Lstp(e/abx)=Lstp(d/acx)$, in the obvious sense (here $e,d$ play the parts of $d,e$ in Claim \ref{vee},
respectively). Now as $d$ extends the type of $e$, $d \ua ab, d \ua ac, d \ud_{ax} bc$. $\{a,b,c\}$ is independent, $U(d/abc) \leq 1$ and as $d \da_{ax} bc$ then 
$U(d/ax) \leq 1$ but $U(e/xab)=U(e/x)=2$ so $e \ua_{ax} b$. So we have $e \ua_{ax} b$,
hence $e \ua_x ab$. But this contradicts our choice of $e$. And so we are done.

\end{proof}

\begin{claim}
\label{2-ind} Let $xRa \land xRb \land xRc$, $Lstp(a/x)=Lstp(b/x)=Lstp(c/x)$, $a,b,c$ 2-independent. Then $a \ua_c b$.
\end{claim}

\begin{proof} Suppose $\{a,b,c\}$ are independent.
Then $b \ud_x c$ and $Lstp(b/x)=Lstp(c/x)$, so we can continue $b,c$ to a Morley sequence.
Let $d,e$ be in this Morley sequence such that $Lstp(d/xa)=Lstp(e/xa)$, and also such that $a \ud_x de$. The first requirement is
achieved by basic counting, and the second since $a \ua_x de \land a \ua_x fg \land a \ua_x hi$, so then as $d,e,f,g,h,i$ are Morley 
over $x$ and $a \ua x$,
$a$ would fork too much for its U-rank 3.  Now $Lstp(d/xb)=Lstp(e/xb)$, as $d,e$ are part of a Morley sequence over $x$ which starts 
with $b$, 
as well as $Lstp(d/xc)=Lstp(e/xc)$. Also $d \ud c \land d \ud b \land d \ud e$ and $e \ud b \land e \ud c \land e \ud d$ as it is an 
indiscernible sequence. So by claim 2, $de \ua b$, $de \ua c$, and $de \ua a$. As $a,b,c$ are independent, $de \ua_a b$ and $de \ua_{ab} c$, so $U(de/abc) \leq 3$, but 
$U(de/x)=4$ so $de \ua_x abc$.  
Notice that all we used to get that $de \ua_x abc$ was $tp(de/xa) \land tp(de/xbc)$. 

We now show we can get such a $d,e$ such that $de \ud_x abc$ and get a contradiction. $de \ud_x a$, 
since chosen as such. As $bcde$ is Morley over $x$, $bc \ud_x de$. Also by rank calculations $bc \ud_x a$, and obviously
$Lstp(de/x)=Lstp(de/x)$, so by the independence theorem over Lascar strong types we get $d^\prime e^\prime$ such that $d^\prime e^\prime
\ud_x bc$ and $tp(d^\prime e^\prime /xbc)=tp(de/xbc)$ and $tp(d^\prime e^\prime /xa)=tp(de/xa)$. So $d^\prime e^\prime$ is exactly as
$de$ with regard to $xa$ and $xbc$, which is what we used in the beginning, and we got the contradiction. \end{proof}

\noindent With this, our Theorem is done. $\square$

\end{proof}

This theorem gives us structure information regarding
IP sequences for a formula which witnesses forking in simple theories. i.e., consequences of the cohabitation of independence property and forking in simple theories. 

\section*{References}

\br

[K]  Byunghan Kim. ``Simplicity, and stability in there.''  The Journal of Symbolic Logic, 66:822-836, 2001. 

\br 
 
[KPi] Byunghan Kim and Anand Pillay. ``Simple theories.'' Annals of Pure and 
Applied Logic, 88:149-164, 1997.

\br

[KPi1] Byunghan Kim and Anand Pillay. ``Around stable forking.'' Fundamenta Mathematicae 170:107-118, 2001

\br
  
[P] Assaf Peretz. \emph{Investigations into the geometry of forking in simple theories}. Ph.D. thesis, UC Berekeley, 2003.

\br
  
  [S1] Saharon Shelah. \emph{Classification Theory}. North-Holland, Amsterdam, The 
Netherlands, 1978.

\br

  [S2] Saharon Shelah. ``Simple unstable theories.'' Annals of Pure and Applied 
Logic, 19:177-203, 1980.

\br

  [W] Frank Wagner. \emph{Simple Theories}. Kluwer Academic Publishers, Dordrecht, 2000.

\end{document}